\documentclass[12pt]{amsart}

\newcommand {\supplus}{\mathop{{\supset}\llap{\raise 
0.5pt\hbox{\normalfont\small+}\hskip 0.5pt}}} 

\newcommand {\subplus}{\mathop{{\subset}\llap{\raise 
0.5pt\hbox{\normalfont\small+}\hskip 0.5pt}}}  

\newcommand {\Cee}    {{\mathbb  C}}

\newcommand {\fb}     {{\mathfrak{b}}}

\newcommand {\fg}     {{\mathfrak{g}}}    %
\newcommand {\fh}     {{\mathfrak{h}}}

\newcommand {\fk}     {{\mathfrak{k}}}

\newcommand {\fn}     {{\mathfrak{n}}}

\newcommand {\fp}    {{\mathfrak{p}}}   %

\newcommand {\fr}     {{\mathfrak{r}}}

\newcommand {\ft}     {{\mathfrak{t}}}

\newcommand {\cal} {\mathcal}

\def \opname#1#2%
  {\expandafter\newcommand \csname #1\endcsname {{\mathop{#2}\nolimits}}}

\newcommand{\rmname}[1]
  {\expandafter\newcommand \csname #1\endcsname {{\operatorname{#1}}}}

\newcommand{\rmnameii}[2]
  {\expandafter\newcommand \csname #1\endcsname {{\operatorname{#2}}}}

\rmname{act}
\rmname{Ad}
\rmname{Add}
\rmname{ad}
\rmname{Alt}
\rmname{alt}
\rmname{Ann}
\rmname{antidiag}
\rmname{Ber}
\rmname{ber}
\rmname{Br}
\rmname{card}
\rmname{ch}
\rmname{Char}
\rmname{cem}
\rmname{cj}
\rmname{Cliff}
\rmname{cntr}
\rmname{codim}
\rmname{Coind}
\rmname{const}
\rmname{col}
\rmname{cork}
\rmname{cpr}
\rmname{diag}
\rmnameii{Div}{div}
\rmname{Def}
\rmname{Der}
\rmname{Dim}
\rmname{End}
\rmname{Even}
\rmname{Ext}
\rmname{gr}
\rmname{Hom}
\rmname{HT}
\rmnameii{Ht}{ht}
\rmname{hwt}
\rmname{Id}
\rmname{id}
\rmname{ind}
\rmname{Ind}
\rmname{Inf}
\rmname{irr}
\rmname{Le}
\rmname{Lie}
\rmname{lwt}
\rmname{mult}
\rmname{Mor}
\rmname{nm}
\rmname{Ob}
\rmname{Odd}
\rmname{Osc}
\rmname{per}
\rmname{Pic}
\rmname{pr}
\rmname{pro}
\rmname{Prime}
\rmname{Proj}
\rmname{prt}
\rmname{pt}
\rmname{Q}
\rmname{qet}
\rmname{qtr}
\rmname{rd}
\rmname{rk}
\rmname{row}
\rmname{Res}
\rmname{salt}
\rmname{Sch}
\rmname{SBr}
\rmname{scalar}
\rmname{Ser}
\rmname{sign}
\rmname{Smbl}
\rmname{spin}
\rmname{ssym}
\rmname{str}
\rmname{st}
\rmname{sgn}
\rmname{sq}
\rmname{symm}
\rmname{supp}
\rmname{Supp}
\rmname{St}
\rmname{Spec}
\rmname{Spm}
\rmname{tr}
\rmname{vpt}
\rmname{weyl}
\rmname{Weyl}
\rmname{Witt}

\opname{vvol}  {{v\hspace{-0.1ex}o\hspace{-0.02ex}l\/}}
\opname{pnt}  {\text{\normalfont pt}}
\opname{Span} {{Span}}
\opname{slim} {\overline{\lim}}
\opname{Vol}  {{V\hspace{-0.55ex}o\hspace{-0.02ex}l\/}}
\opname{Par}  {{P\hspace{-0.3ex}a\hspace{-0.05ex}r\/}}

\rmname{Mat}
\rmname{Bil}
\rmname{Diff}
\rmname{Ker}
\rmname{Herm}
\rmname{Coker}
\rmname{Conn}
\rmname{Covect}
\rmname{Vect}
\rmname{Int}

\rmnameii {IM} {Im}
\rmnameii {RE} {Re}

\opname{Aut} {{A\hspace{-0.2ex}u\hspace{-0.1ex}t\/}}
\opname{GL} {{G\hspace{-0.3ex}L}}
\opname{SL} {{S\hspace{-0.3ex}L}}
\opname{Exp} {{E\hspace{-0.2ex}x\hspace{-0.1ex}p\/}}
\opname{GQ} {{G\hspace{-0.2ex}Q}}
\opname{OSp} {{O\hspace{-0.25ex}S\hspace{-0.15ex}p\/}}
\opname{Out} {{O\hspace{-0.25ex}u\hspace{-0.15ex}t\/}}
\opname{Spp} {{S\hspace{-0.2ex}p\/}}
\opname{SpO} {{S\hspace{-0.2ex}p\hspace{-0.02ex}O\/}}
\opname{Pe} {{P\hspace{-0.25ex}e\/}}
\opname{SPe} {{S\hspace{-0.25ex}P\hspace{-0.25ex}e\/}}
\opname{Spin} {{S\hspace{-0.25ex}p\hspace{-0.05ex}i\hspace{-0.1ex}n\/}}
\opname{Iso} {{I\hspace{-0.25ex}s\hspace{-0.1ex}o\/}}
\opname{SSPe} {{S\hspace{-0.25ex}S\hspace{-0.15ex}P\hspace{-0.25ex}e\/}}
\opname{PeU} {{P\hspace{-0.25ex}e\hspace{-0.1ex}U\/}}
\opname{QU} {{Q\hspace{-0.15ex}U\/}}
\opname{U} {{U\/}}

\opname{cGQ} {{\cal G \hspace{-0.2em} Q \/}}
\opname{cSL} {{\cal S \hspace{-0.2em} L \/}}
\opname{cGL} {{\cal G \hspace{-0.2em} L \/}}
\opname{cOSp} {{\cal O \hspace{-0.2em} S \hspace{-0.3em} \it p\/}}
\opname{cPe} {{\cal P \hspace{-1.5pt} \it e\/}}
\opname{cVect} {{\cal V \hspace{-1.5pt} \it e\hspace{-0.1ex}c\hspace{-0.1ex}t\/
}}
\opname{cVol} {{\cal V \hspace{-1.5pt} \it o\hspace{-0.1ex}l\/}}
\opname{cAut} {{\cal A \hspace{-0.2em} \it u\hspace{-0.1em}t\/}}
\opname{cCovect} {{\cal C \hspace{-1.5pt}
     \it o\hspace{-0.1ex}v\hspace{-0.1ex}e\hspace{-0.1ex}c\hspace{-0.1ex}t\/}}
\opname{CW} {{C\hspace{-0.15ex}W}}

\newcommand {\ev} {{\bar0}}
\newcommand {\od} {{\bar1}}

\newcommand {\tto} {\longrightarrow}

\newcommand {\bcdot}   {\mathbin{\hbox{\raise.4ex\hbox{\bf.}}}} 

\newcommand {\secno} {}
\newcommand {\ssecfont} {\normalfont\bf}

\newtheorem{Theorem}{\secno Theorem}

\newtheorem{Lemma}[Theorem]{\secno Lemma}

\newtheorem{Corollary}[Theorem]{\secno Corollary}

\newenvironment {th*}[1]
    {\gdef\thname{#1} \begin{thn}}%
    {\end{thn}}
\newtheorem{thn}[Theorem] {\thname}

\theoremstyle{definition}

\newenvironment {ex*}[1]
    {\gdef\thname{#1} \begin{exn}}%
    {\end{exn}}
\newtheorem{exn}[Theorem]{\thname}

\theoremstyle{remark}

\newenvironment {rem*}[1]
    {\gdef\thname{#1} \begin{remn}}%
    {\end{remn}}
\newtheorem{remn}[Theorem]{\thname}

\newcommand {\ssec}{\subsection*}

\newcommand {\ssbegin}[2]
  {\def \secno {\gdef \secno {}{\ssecfont #1. }}
   \begin{#2}}

\begin{document}

\title[Irreducible representations of solvable Lie 
superalgebras]{Irreducible representations of solvable Lie 
superalgebras}

\author{Alexander Sergeev}

\address{Dept. of Math., Univ. of Stockholm, Roslagsv. 101, 
Kr\"aftriket hus 6, S-106 91, Stockholm, Sweden; 
mleites@matematik.su.se (on leave of absence 
from Balakovo Institute of Technichnics, Technology and Control)}

\thanks{I am thankful to I.~Shchepochkina and D.~Leites for help.}

\begin{abstract} The description of irreducible finite dimensional 
representations of finite dimensional solvable Lie superalgebras over 
complex numbers given by V.~Kac is refined.  In reality these 
representations are not just induced from a polarization but twisted, 
as infinite dimensional representations of solvable Lie algebras.  
Various cases of irreducibility (general and of type Q) are classified.
\end{abstract}



\maketitle

\section*{Introduction}

Hereafter the ground field is $\Cee$ and all
the modules and superalgebras are finite dimensional; $\fg$ is a
solvable Lie superalgebra.

A description of irreducible representations of solvable Lie 
superalgebras given in \cite{K} (Theorem 7) contains an error.  In 
reality to give such a description one has to imitate the description 
of infinite dimensional solvable Lie algebras \cite{D}, i.e., we must 
consider {\it twisted} induced representations.  In what follows I 
give a correct description of irreducible representations of solvable 
Lie superalgebras.  I also show where a mistake crept into \cite{K} 
and give a counterexample to Theorem 7 from \cite{K}.

The proof given in what follows was delivered at Leites' {\it Seminar 
on Supermanifolds} in 1983 and is preprinted in \cite{L} in a form 
considerably edited by I.~Shchepochkina and D.~Leites.  My 
acknowledgements are due to them and also to the Department of 
Mathematics of Stockholm University that financed publication of 
\cite{L}.

\section*{\S 1. Main result}

\ssec{1.1. Polarizations} Set
$$
L=\{\lambda \in \fg ^{*}:\lambda (\fg _{\od})=0\text{ and }
\lambda ([\fg _{\ev}, 
\fg _{\ev}])=0\}. 
$$
It is convenient to represent the functionals from $L$ as elements of an 
isomorphic space
$$
L_{0}=(\fg _{\ev}/[\fg _{\ev}, \fg _{\ev}])^{*}.
$$ 
The natural isomorphism of these spaces is given by the formula
$$
' : L_{0}\tto L,\; \; \lambda ' \mid \fg _{\ev}=\lambda , 
\; \; \lambda ' \mid \fg _{\od}=0, \eqno{(')}
$$
where $\lambda$ denotes a character and also, by the usual abuse of 
language, the (1, 0)-dimensional representation of the Lie algebra $\fg_{\ev}$
determined by the character $\lambda$.  Every functional $\lambda \in L$
determines a symmetric  form
$f_{\lambda}$ on $\fg _{\od}$ by the formula
$$
f_{\lambda}(\xi _{1}, \xi _{2})=\lambda ([\xi _{1}, \xi _{2}]).
$$
A subalgebra $\fh \subset \fg $ is called a {\it 
polarization} for $\lambda \in L$ (and is often denoted by $\fp(\lambda)$) if 
$$
\lambda ([\fh , \fh])=0,\; \; \fh \supset \fg _{\ev}\; \; \text{and}\; \fh
_{\od}\; \text{ is a maximal isotropic subspace for}\; f_{\lambda}. 
$$

\begin{Lemma} For every $\lambda\in L$ there exists a polarization $\fh$. 
\end{Lemma}

If $\fh $ is a polarization for $\lambda \in L$, then, 
clearly, $\lambda $ determines a (1, 0)-dimensional representation of
$\fh$.

\ssec{1.2. Twisted representations} Let $\fh \subset \fg $ be a Lie
subsuperalgebra that contains 
$\fg _{\ev}$. Define a functional $\theta _{\fh}\in L$ by setting
$$
\theta _{\fh (g)}=\left \{\begin{matrix}
\frac{1}{2}\str_{\fg /\fh}(\ad~g)&\text{for}\,
g \in \fg _{\ev}\cr
0&\text{for}\, g\in \fg _{\od} \end{matrix}\right .
$$
Note that $\theta _{\fh}([\fh , \fh])=0$ by definition of the
supertrace. Therefore, $\theta _{\fh}$ is a character of a $(1, 0)$-dimensional
representation of $\fh $.

Let $\fh $ be a polarization for $\lambda \in L$. Define the twisted (by
the character $\theta _{\fh}$) induced and coinduced representations by setting
$$
\begin{gathered}
I^{\fg}_{\fh}(\lambda)=\Ind^{\fg}_{\fh}(\lambda
+\theta _{\fh})=U(\fg)\otimes _{U(\fh)}(\lambda +\theta
_{\fh});\\
CI^{\fg}_{\fh}(\lambda)=\Coind^{\fg}_{\fh}(\lambda -\theta _{\fh})=
\text{Hom}_{U(\fh)}(U(\fg), \lambda -\theta _{\fh}).
\end{gathered}
$$

\begin{Lemma} {\em 1)} $I^{\fg}_{\fh}(\lambda)$ is finite dimensional 
and irreducible.

{\em 2)} $I^{\fg}_{\fh}(\lambda)$ does not depend on the choice of a 
polarization $\fh $; therefore, notation $I(\lambda)$ 
($=I^{\fg}_{\fh}(\lambda)$ for some $\fh$) is well-defined.

{\em 3)} $CI(\lambda)\cong I(\lambda )$.
\end{Lemma}

\ssec{1.3. Main Theorem} Let $Z=\{(\lambda , \fh): \lambda \in L$ and let $\fh 
$
be a polarization for $\lambda \}$. Define an equivalence relation on $Z$
by setting
$$
(\lambda , \fh)\sim (\mu , \ft)\Longleftrightarrow \lambda -\theta _{\fh}=
\mu -\theta _{\ft}
$$
Clearly, this relation is well-defined.

Recall (\cite{BL}) that the representation of a Lie superalgebra $\fg$ 
is called irreducible of {\it $G$-type} if it has no invariant 
subspaces; it is called irreducible of {\it $Q$-type} if it has no 
invariant sub{\it super}spaces.

\begin{Theorem} {\em 1)} Every irreducible finite dimensional 
representation of $\fg $ is isomorphic up to application of the change 
of parity functor $\Pi $ to a representation of the form $I(\lambda)$ 
for some $\lambda $.

{\em 2)} The map $\lambda \tto I(\lambda)$ is (up to $\Pi$) a 1-1
correspondence
between elements of $L$ and the irreducible finite dimensional
representations of $\fg $.

{\em 3)} Let $(\lambda , \fh), (\mu , \ft)\in Z$. Then 
$\Ind^{\fg}_{\fh}(\lambda)\cong \Ind^{\fg}_{\ft}(\mu)$ if and
only if $(\lambda , \fh)\sim (\mu , \ft)$.

{\em 4)} If $\rk f_{\lambda} $ is even, then $I(\lambda)$ is a 
$G$-type representation; if 
$\rk f_{\lambda}$ is odd, then $I(\lambda)$ is a $Q$-type 
representation.
\end{Theorem}

\section*{ \S 2. Prerequisites for the proof of Main theorem}

Let $\fk\subset \fg $ be a subsuperalgebra, $\codim~\fk =\varepsilon$, 
and $\mu $ the character of the representation of $\fg _{\ev}$ in $\fg 
/\fk$.  For the definition of the isomorphism $'$ see sec. 1.1.

\ssbegin{2.1}{Lemma} $\mu'$ is a character of $\fg$.
\end{Lemma}

\begin{proof} Let $\xi \in \fg$ and $\xi \not\in \fk$.  Since in 
$\fg/\fk$ there is a $\fk$-action that coincides with $\mu '|\fk$, 
it suffices to prove that $\mu ' ([\fk, \xi])=\mu '([\xi , \xi])=0$.  
By the Jacobi identity $[[\xi , \xi], \xi]=0$ which proves that $\mu ' 
([\xi , \xi])=0$.  Since $p(\mu)=\ev$, we have $\mu ([\xi , \xi])=0$.  
Let $\eta \in \fk_{\od}$.  Then $[[\eta , \xi], \xi]=\frac{1}{2}[\eta, [\xi ,
\xi]]\in \fk$ and, therefore, $\mu ([\fk _{\od}, \xi])=0$.
\end{proof}

\ssbegin{2.2}{Corollary} Let $\fg _{\ev}\subset \fh \subset 
\fk\subset \fg $ so that $\dim \fk/\fh =\varepsilon$.  Let $\lambda $ 
be the weight of a vector in the $\fg _{\ev}$-module $\fg /\fh $.  
Then $\lambda '$ is a character of $\fk $. 
\end{Corollary}

\ssbegin{2.3}{Corollary} {\em (\cite{K}).} A Lie superalgebra 
$\fg =\fg _{\ev}\oplus \fg _{\od}$ is solvable if and
only if so is $\fg _{\ev}$.
\end{Corollary}

\ssbegin{2.4}{Corollary} If $\fh $ is a polarization for $\lambda$, 
then $\fh $ is also polarization for $\lambda +\alpha \theta _{\fh}$ 
for any $\alpha \in \Cee^{*}$.
\end{Corollary}

Let us recall three well-known lemmas.

\ssbegin{2.5}{Lemma} {\em (see \cite{K}).} Let $\fg =\fk \oplus \Span(g)$, 
where $\fk$ is a Lie subalgebra and $p(g)=\od$.  If $(V, \rho)$ is an 
irreducible representation of $\fk$ in a superspace $V$, then 
$W=\Ind^{\fg}_{\fk}(V)$ is reducible if and only if $V$ admits a 
$\fg$-module structure that extends $\rho$.
\end{Lemma}

\ssbegin{2.6}{Lemma} {\em (see \cite{K}).} Let $\fk\subset \fg$ be a Lie 
subsuperalgebra, $\dim \fg /\fk =(0, 1)$.  If $W$ is an irreducible 
$\fg$-module and $V\subset W$ is an irreducible proper 
$\fk$-submodule, then $W=\Ind^{\fg}_{\fk}(V)$.
\end{Lemma}

\ssbegin{2.7}{Lemma} Let $W$ be a finite dimensional $\fg $-module, 
$f$ a symmetric $\fg $-invariant form on $W$ and $V$ a $\fg 
$-invariant isotropic subspace.  Then there exists a maximal $\fg 
$-invariant $f$-isotropic subspace in $W$ containing $V$.
\end{Lemma}

{\bf Proof} follows from linear algebra. \qed

\ssbegin{2.8}{Corollary} Lemma $1.1$ holds. \qed
\end{Corollary}

\section*{\S 3. Description of irreducible modules}

\ssbegin{3.1}{Proposition} Let $\lambda \in L$, let $\fh =\fg_{\ev}\oplus 
\fp$ be a polarization for $\lambda$, $\fn$ the kernel of $f_{\lambda}$ 
and $F\subset \fp$ a subspace such that $\fp=F\oplus \fn, \xi _{0}\in 
\fp^{\bot}$ and $\xi _{0}\not\in \fp $ (if $\rk~f_{\lambda}$ is 
even, then we set $\xi _{0}=0)$, $v$ a generator of 
$\Ind^{\fg}_{\fk}(\lambda)$.  If $u\in 
\Ind^{\fg}_{\fk}(\lambda)$ and $Fu=0$, then $u\in \Span(v, \xi_{0}v)$.
\end{Proposition}

{\bf Proof} is carried out by induction on $\rk ~f_{\lambda}$.  If 
$\rk ~f_{\lambda}=0$, then $F=0$ and the statement is obvious.

Let $\rk ~f_{\lambda}>0$.  Select a subalgebra $\fh\supset\fb$ such 
that and $\dim\fg_{\od}/\fh _{\od}=1$.  The two cases are possible:
$\fh_{\od}^\perp\not\subset\fh_{\od}$ and 
$\fh_{\od}^\perp\subset\fh_{\od}$.

\underline{i) $\fh_{\od}^\perp\not\subset\fh_{\od}$}.  Then $\fg 
_{\od}=\fh _{\od}\oplus\Span(\xi)$, where $\xi\perp \fh_{\od}$.  
Hence, $\xi\perp \fb_{\od}$ and $\xi\not\in \fb_{\od}$.  Therefore, we 
may assume that $\xi=\xi_{0}$ and $\rk ~f_{\lambda}$ is an odd number.  
Clearly, $\fb$ is a polarization for the restriction $f_{\lambda}'$ 
onto $\fh_{\od}$ and $\rk ~f_{\lambda}'$ is an even number.  Further on,
$$
\Ind^{\fg}_{\fb}(\lambda)=\Ind^{\fh}_{\fb}(\lambda)\oplus 
\xi _{0}\Ind^{\fg}_{\fb}(\lambda).
$$

Let $u=u_{0}+\xi _{0}u_{1}\in \Ind^{\fg}_{\fb}(\lambda)$ and $pu=0$ 
for any $p\in F$.  Then
$$
0=pu=pu_{0}+[p, \xi _{0}]u_{1}-\xi _{0}pu_{1}, 
$$
therefore, $pu_{1}=0$.  By induction, $u_{1}\in \Span(v)$, where $v$ 
is the generator of $\Ind^{\fg}_{\fb}(\lambda)$.  Since $\xi 
_{0}\perp\fb_{\od}$, it follows that $[p, \xi_{0}]u_{1}= 
f_{\lambda}(p, \xi _{0})u_{1}=0$.  Therefore, $pu_{0}=0$ and $u_{0}\in 
\Span(v)$. Hence, $u\in \Span(v, \xi _{0}v)$.

Let us show now that the weight of $\xi _{0}v$ with respect to 
$\fg_{\ev}$ is also equal to $\lambda$. Indeed, since $[x, 
\xi_{0}]\perp \fb_{\od}$ for any $x\in \fg_{\ev}$, it follows that 
$[x, \xi_{0}]=\mu(x)\xi_{0}+p$ for some $p\in \fb_{\od}$. Furthermore, 
$[x, [\xi_{0}, \xi_{0}]=2[[x, \xi_{0}], \xi_{0}]$; hence, 
$$
\begin{array}{l}
0=\lambda([x, [\xi_{0}, \xi_{0}])=2\lambda([\mu(x)\xi_{0}+p, 
\xi_{0}])=\\
2\mu(x)\lambda([\xi_{0}, \xi_{0}])+2\lambda([p, \xi_{0}])=2\mu(x).\\
\end{array}
$$
So, $\mu(x)=0$ and the weight of $\xi _{0}v$ is equal to $\lambda$. 

\underline{ii) $\fh_{\od}^\perp\subset\fh_{\od}$}.  Then the {\it 
restriction} of the form $f_{\lambda}$ onto $\fh_{\od}$ is of rank by 
2 less than that of $f_{\lambda}$ itself.

Select $\xi\perp\fh_{\od}$, $\eta\in \fh_{\od}^\perp$ and 
$F_{1}\subset\fb_{\od}$ so that
$$
F=F_{1}\oplus\Span(\eta); \; \xi\perp F_{1}; \; f_{\lambda}(\xi , 
\eta)\neq 0.
$$ 
Let 
$$
u=u_{0}+\xi u_{1}, \text{ where } u\in \Ind^\fg_{\fb}(\lambda), \; 
u_{0}, u_{1}\in \Ind^\fh_{\fb}(\lambda)\text{ and } pu=0 \text{ for 
any }p\in F. 
$$  
Then
$$
0=pu=pu_{0}+[p, \xi]u_{1}+\xi pu_{1},
$$
hence, $pu_{1}=0$ and by induction $u_{1}\in \Span(v, \xi _{0}v)$. 
Thanks to i) $[p, \xi]u_{1}=f_{\lambda}(p, \xi)u_{1}$ and if $p\in 
F_{1}$, then $f_{\lambda}(p, \xi)u_{1}=0$; hence, $pu_{0}=0$ for any $p\in 
F_{1}$. By induction we deduce that $u_{0}\in \Span(v, \xi _{0}v)$.
Further, 
$$
0=\eta u=\eta u_{0}+\eta \xi u_{1}=[\eta, \xi]u_{1}=f_{\lambda}(\eta, 
\xi)u_{1}
$$ 
and since $f_{\lambda}(\eta, \xi)\neq 0$, then $u_{1}=0$ and 
$u=u_{0}\in\Span(v, \xi_{0}v)$.  \qed

\ssbegin{3.2}{Corollary} If $\fh =\fg _{\ev}\oplus \fp$ is a 
polarization for $\lambda$, then $\Ind^{\fg}_{\fh}(\lambda)$ 
is an irreducible module.
\end{Corollary}

{\bf Proof} Irreducibility is equivalent to the absence of singular 
vectors that do not lie in $\Span(v, \xi _{0}v))$. \qed

\ssbegin{3.3}{Corollary} heading 1) of Lemma 1.2 holds.
\end{Corollary}

{\bf Proof} follows from the definition of $I^{\fg}_{\fh}(\lambda)$
and sect. 2.4. \qed

\ssbegin{3.4}{Corollary} Let $U$ be an irreducible finite dimensional 
$\fg $-module.  Then $U=\Ind^{\fg}_{\fh}(\lambda)$ for some $\lambda 
\in L$ and a polarization $\fh $.
\end{Corollary}

{\bf Proof} will be carried out by induction on $\dim \fg _{\od}$.  If 
$\fg =\fg _{\ev}$, then this is Lie's theorem.  Let $\fk \subset \fg $ 
and $\dim \fg _{\od}/\fk _{\od}=1$.

Let $U$ be irreducible as a $ \fk $-module.  Then there exist $\lambda 
\in L$ and a polarization $\fh \subset \fk$ for $\lambda \in L$ such 
that $U=\Ind^{\fk}_{\fh}(\lambda)$.  If $\fh $ were a polarization for 
$\lambda $ in $\fg $, too, then by Corollary 3.2 the representation
$$
\Ind^{\fg}_{\fh} (\lambda)=
\Ind^{\fg}_{\fk} (\Ind^{\fk}_{\fh} (\lambda))
$$
would have been irreducible contradicting Lemma 2.5.

Let $\hat{\fh}\supset \fh $ be a polarization for $\lambda $ in $\fg $ 
and $\xi \in \hat{\fh}$ so that $\xi \not\in \fh $.  If $v$ is a 
generator of $\Ind^{\fk}_{\fh}(\lambda)$ and $p\in\fh_{\od}$, then
$$
\begin{gathered}
p\xi v=[p, \xi]v=f_{\lambda}(p, \xi)v=0,\\
\xi \xi v=\frac12[\xi , \xi]v=\frac12f_{\lambda}(\xi , \xi)v=0.
\end{gathered}
$$
Therefore, there exists a non-zero $\fg $-module homomorphism 
$\Ind^{\fg}_{ \hat{\fh}}\lambda \tto \Ind^{\hat{\fk 
}}_{\fh}(\lambda)=U$ and since both modules are irreducible, this is 
an \lq\lq odd isomorphism", i.e., the composition of an isomorpism 
with the change of parity.

Now let $U$ be reducible as a $\fk$-module.  Then by Lemma 2.6 
$U=\Ind^{\fg}_{\fk}V$, and, by induction, 
$V=\Ind^{\fk}_{\fh}(\lambda)$ for a polarization $\fh \subset \fk $ 
and $\lambda \in L $.  If $\fh $ is not a polarization for $\lambda $ 
in $\fg$, then let $\hat{\fh}\supset \fh $ be a polarization.  We have 
a non-zero $\fg $-module homomorphism $U=\Ind^{\fg}_{\fh}(\lambda) 
\tto \Ind^{\fg}_{\hat{\fh}}(\lambda)$ and since both modules are 
irreducible, this is an isomorphism which is impossible because $\dim 
\Ind^{\fg }_{\hat{\fh}}(\lambda)<\dim \Ind^{\fg}_{\fh}(\lambda)$.  
Therefore, $\fh $ is a polarization for $\lambda $ in $\fg $ and 
$U=\Ind^{\fg}_{\fh}(\lambda)$.  \qed

\ssbegin{3.5}{Corollary} Heading $1)$ of Theorem holds.
\end{Corollary}

\ssec{3.6. A subsuperalgebra subordinate for $l\in L$}
Recall, see [K], that if
$$
\fg _{l}=\{g\in \fg \mid l([g, g_{1}])=0\; \text{ for all } \; \; 
g_{1}\in \fg \},
$$
then a subalgebra $\fp \subset \fg $ is called {\it subordinate to } 
$\lambda $ if $l([\fp , \fp])=0$ and $\fp \supset \fg _{l}$.

\begin{Corollary} Let $l\in L$, $\fb$ a subalgebra subordinate 
to $l$.  Then $\Ind^{\fg}_{\fb}(l)$ is 
irreducible if and only if $\fb$ is a polarization for $l$.
\end{Corollary}

\section*{\S 4. Classification of modules $\Ind^{\fg}_{\fh}(\lambda)$}

\ssbegin{4.1}{Lemma} If $(\lambda , \fh)\sim (\mu , \ft)$, then $\fh $ 
is a polarization for $\mu $.
\end{Lemma}

\begin{proof} By 2.4 $\fh $ is a polarization for $\lambda -\theta 
_{\fh}$.  Since $\lambda -\theta _{\fh}=\mu -\theta _{\ft}$, then $\fh 
$ is a polarization for $\mu -\theta _{\ft}$ and since $\ft $ is a 
polarization for $\mu -\theta _{\ft}$, it follows that $\dim \fh =\dim 
\ft $.  Let $\fn$ be the kernel of $f_{\mu -\theta _{\ft}}$.  Then 
$\fh\cap\ft \supset \fg _{\ev}\oplus \fn$ and, therefore, $\fg /\fb $ 
is a quotient of $\fg _{\od}/\fn $; hence by sect.  2.2 we have 
$\theta _{\ft}([\fh, \fh])=0$.

Therefore, $\fh $ is completely isotropic with respect to $f_{\mu}$ and
since $\dim \fh =\dim \ft $, we see that $\fh $ is a polarization for 
$\mu $. 
\end{proof}

\ssec{4.2.  Proof of heading 3) of Theorem} Let $(\lambda , \fh) \sim 
(\mu , \ft)$.  We will carry the proof out by induction on $k=\dim 
\fh /(\fh \cap \ft)$.  If $k=0$ the statement is obvious.  Let 
$k=1, $ then, obviously, $\dim \ft /(\fh \cap \ft)=1$.  Consider 
the space $\fh +\ft $.  By Lemma 4.1 $\ft $ is a polarization for 
$\lambda $ and, therefore, the kernel of $f_{\lambda}$ on the 
subspace $\fh +\ft $ is equal to $\fh \cap \ft $.  Let $\xi \in \fh 
$ and $\eta \in \ft $ be such that $\bar{\xi}\in \fh /(\fh \cap 
\ft), \bar{\xi}\neq 0$ and $\bar{\eta}\in \ft /(\fh \cap \ft), 
\bar{\eta}\neq 0$.  We may assume that $f_{\lambda}(\xi , \eta)=1$.

Let $v$ be a generator of $\Ind^{\fg}_{\fh}(\lambda)$.  Then for $r\in 
\fh \cap \ft $ we have
$$
r\eta v=[r, \eta]v=\lambda ([r, \eta])v=0,\; \; \; 
\eta \eta v=\frac12[\eta , \eta]v=\frac12\lambda ([\eta
, \eta])v=0,
$$
i.e., $\ft _{\od}(\eta v)=0$ and, therefore, there exists a non-zero 
homomorphism $\Ind^{\fg}_{\ft}(\mu')\tto \Ind^{\fg}_{\fh}(\lambda)$, 
where $\mu'$ is the weight of $\eta v$.

Since $\Ind^{\fg}_{\fh}(\lambda)$ is irreducible and $\dim 
\Ind^{\fg}_{\ft}(\mu')=\dim \Ind^{\fg}_{\fh}(\lambda)$, this 
homomorphism is an isomorphism.  Let $g\in \fg _{\ev}$.  Then
$$
g(\eta v)=\eta (gv)+[g, \eta]v=[\lambda +\tr_{\ft 
/(\ft\cap\fh)}\ad~g]\eta v,
$$
i.e., $\mu ' =\lambda -\str _{\ft /(\ft \cap \fh)}\ad g$.

Since $\lambda \in L $, it follows that
$$
\begin{gathered}
0=\lambda ([g, [\xi , \eta]])=\lambda ([[g, \xi], \eta])+\lambda 
([\xi, [g, \eta ]])=\\
=-(\str _{\ft /\ft \cap \fh}\ad~g+\str _{\fh /\ft 
\cap \fh}\ad~g)\lambda ([\xi , \eta]).
\end{gathered}
$$
Since $\lambda ([\xi , \eta]=1$, it follows that
$$
\str _{\ft /\ft\cap \fh}\ad~g=
-\str _{\fh/\ft \cap \fh}\ad~g, 
$$
and
$\mu' =\lambda -\theta _{\fb}-\theta _{\ft}$, i.e., 
$\Ind^{\fg}_{\ft}(\mu)\cong \pi (\Ind^{\fg}_{\ft}(\lambda))$.

Let $k>1$.  Let $\fh =\fg _{\ev}+P$ and $\ft =\fg _{\ev}\oplus Q$ such 
that $P\cap Q\subset F\subset P$, $F\neq P$ and $F\neq P\cap Q$, where 
$F$ is a $\fg_{\ev}$-submodule in $\fg _{\od}$.  Set 
$R=F+(F^{\bot}\cap Q)$.  It is not difficult to verify that $\fr=\fg 
_{\ev}\oplus R$ is a polarization for $\lambda$.  Set
$$
\nu (x)=\lambda (x)+\str_{P/(P\cap R)}(\ad~x).
$$
Since $P/(P\cap R)$ is a subquotient in $\fg _{\od}/\fn $, it follows 
that $\nu ([R, R])=0$.  The same arguments as in Lemma 4.1 show that 
$\fr$ is a polarization for $\nu$.

Since $P\cap R\supset F\supset P\cap Q$, then
$\dim R/(P\cap R)<\dim P/(P\cap Q)$.

Further, the diagram of inclusions
$$
\begin{matrix}
P~\cap ~R& \tto &~P\cr
\Big \downarrow &&\Big \downarrow \cr
 R&\tto &\fg _{\od}\end{matrix}
$$ 
shows that
$$
2\theta_{\fh} (x)+\str _{P/(P~\cap~R)}\ad~(x)=2\theta_{\fr}(x)+\str 
_{R/(R~\cap ~P)}\ad~ (x).
$$
By duality, there exists a pairing 
$$
(P/(P\cap R))\times (R/(P\cap R))\tto \Cee
$$
and since $\str _{P/(P~\cap ~R)}\ad(x)=-\str _{R/(P~\cap ~R)}\ad 
~(x)$, then $\str _{R/(P~\cap ~R)}\ad(x)=\theta _{\fh}(x)-\theta 
_{\fr}(x)$.

Thus,
$$
\nu (x)-\theta _{\fr}(x)=\lambda (x)-\str _{R/(P~\cap 
~R)}\ad(x)-\theta _{R}(x)=\lambda (x)-\theta _{\fh}(x),
$$
i.e., $(\lambda , \fh)\sim (\nu , \ft)$ and, by induction,
$\Ind^{\fg}_{\fh}(\lambda)=\Ind^{\fg}_{\fr}(\nu)$. Besides, $\nu
-\theta_{\fr}=\lambda -\theta _{\fh}=\mu -\theta_{\ft}$ and
$Q\cap R\supset F^{\bot} \cap Q<\supset Q\cap P$; therefore,
$$
\dim Q/Q\cap R\leq \dim \bar{Q}=\dim Q/(Q\cap P).
$$
By induction, $\Ind^{\fg}_{\fr}(\nu)\cong 
\Ind^{\fg}_{\ft}(\mu)$, therefore,
$\Ind^{\fg}_{\fh}(\lambda)=\Ind^{\fg}_{\ft}(\mu)$.

Conversely, let $\Ind^{\fg}_{\fh}(\lambda)\cong\Ind^{\fg}_{\ft}(\mu)$.  
Then $\lambda =\mu +\alpha _{1}+\ldots +\alpha _{k}$, where the 
$\alpha _{i}$ are the weights of $\fg _{\od}/Q$.  Therefore, by sect.  
2.2 $\lambda ([Q, Q])=0$ and since $\dim 
\Ind^{\fg}_{\fh}(\lambda)=\dim \Ind^{\fg}_{\ft}(\mu)$, it folows that 
$\dim \fh =\dim \ft $ and $\fh $ is a polarization for $\lambda$ and 
for $\mu ' =\lambda -\theta _{\fh}+\theta _{\ft}$, too.  Since $\mu ' 
-\theta _{\ft}=\lambda -\theta _{\fh}$, then by the above 
$\Ind^{\fg}_{\ft}(\mu')=\Ind^{\fg}_{\fh}(\lambda)= 
\Ind^{\fg}_{\ft}(\mu)$.  Let $v'\in\Ind^{\fg}_{\ft}(\mu ')$ be a 
generator of $\Ind^{\fg}_{\ft}(\mu ')$ and $Qv=0$.  By 4.1 $v' \in 
\Span(v, \xi _{0}v)$, therefore, $\mu ' =\mu $ and $(\mu , \ft) \sim 
(\lambda , \fh)$.  \qed

\ssbegin{4.3}{Corollary} Heading $2)$ of Theorem and heading $2)$ of 
Lemma 1.2 hold.
\end{Corollary}

\begin{proof} Due to sect. 2.4 it is clear that $\fh $ is a polarization for 
$\lambda +\theta _{\fh}$ and, therefore, $I^{\fg}_{\fh}(\lambda)$ is 
irreducible.  If $\ft $ is another polarization for $\lambda$, then by 
sect. 4.2
$$
I^{\fg}_{\ft}(\lambda)=\Ind^{\fg}_{\ft}(\lambda
+\theta _{\ft})=\Ind^{\fg}_{\fh}(\lambda +\theta
_{\fh})=I^{\fg}_{\fh}(\lambda).
$$
If $U$ is irreducible, then by sect.  3.4 $U\cong 
V^{\fg}_{\fh}(\lambda)$ for some $\lambda $ and $\fh $.

If $I(\lambda)=I(\mu)$, then $\Ind^{\fg}_{\fh}(\lambda +\theta
_{\fh})\cong \Ind^{\fg}_{\ft}(\mu +\theta
_{\ft})$ and by sect. 4.2 
$$
\lambda =\lambda +\theta _{\fh}-\theta _{\fb}=
\mu +\theta _{\ft}-\theta _{\ft}=\mu . \qed
$$

{\bf Proof of heading 3) of Lemma 1.2} Let us prove now that 
$I^\fg_{\fb}(\lambda)\cong CI^\fg_{\fb}(\lambda)$.  To this end, make 
use of the isomorphisms $(I^\fg_{\fb}(\lambda))^*\cong 
CI^\fg_{\fb}(-\lambda)$ and $(I^\fg_{\fb}(\lambda))^*\cong 
I^\fg_{\fb}(-\lambda+2\theta_{\fb})$. The first of these isomorphisms 
follows from the definitions of the induced and coinduced modules.

Let us prove the other one. Select a basis $\xi_{1}$,  \dots , 
$\xi_{n}$ in the complement to $\fb_{\od}$ in $\fg_{\od}$ and 
consider the following filtration of $I^\fg_{\fb}(\lambda)$:
$$
I_{0}=\Span(v)\subset I_{1}=\Span(v; \xi_{1}v,  \dots , 
\xi_{n}v)\subset \dots \subset I_{n}=I^\fg_{\fb}(\lambda).
$$
It is clear that the elements $\xi$ can be chosen so that each 
$I_{k}$ is a $\fg_{\ev}$-module. Let $l\in (I^\fg_{\fb}(\lambda))^*$ 
be such that $l(I_{n})\neq 0$ while $l(I_{n-1})=0$. Then it is easy 
to verify that $\fb_{\od}l=0$ and the weight $l$ with respect to 
$\fg_{\ev}$ is equal to $-\lambda+2\theta_{\fb}$. Therefore, there 
exists a nonzero homomorphism $\varphi: 
I^\fg_{\fb}(-\lambda+2\theta_{\fb})\tto (I^\fg_{\fb}(\lambda))^*$. 

Since the dimensions of these modules are equal and the first of them 
is irreducible, $\varphi$ is an isomorphism. Hence,
$$
CI(\lambda)= CI^\fg_{\fb}(\lambda-\theta_{\fb})\cong 
(I^\fg_{\fb}(-\lambda+\theta_{\fb}))^*=I^\fg_{\fb}(\lambda-
\theta_{\fb}+2\theta_{\fb})= I(\lambda).
$$
\end{proof}

\section*{\protect \S 5. An example}

Set $\fg _{\ev}=\Span(x, y, z, u)$, $\fg _{1}=\Span(\xi _{-2}, \xi
_{-1}, \xi _{1}, \xi _{2})$ and
let the nonzero brackets be:
$$
\begin{gathered}
{}[\xi _{-1}, \xi _{1}]=u, \quad [\xi _{-2}, \xi _{2}]=u, \quad [\xi _{1}, \xi
_{2}]=y\\
{}[\xi _{-1}, \xi _{2}]=z, \quad [\xi _{2}, \xi _{2}]=x,\\
[y, \xi _{-2}]=-\frac{1}{2}\xi _{-2}, \quad [y, \xi _{2}]=\frac{1}{2}\xi 
_{1}\\
{}[z, \xi _{1}]=\frac{1}{2}\xi _{-2}, \quad [z, \xi _{2}]=-\frac{1}{2}\xi 
_{1}\\
{}[x, y]=-y, \quad [x, z]=z, \quad [x, \xi _{-1}]=\xi _{-1}, 
\quad [x, \xi _{1}]=-\xi_{1}.
\end{gathered}
$$

It is not difficult to verify that $\fg =\fg _{\ev}\oplus \fg _{\od}$ 
is a Lie superalgebra.  It is solvable since so is $\fg _{\ev}$.  
Consider $\lambda =\lambda _{1}u^*+\lambda _{2}x^*$ with $\lambda 
_{1}\neq 0$.  Then $\fh =\fg _{\ev}\oplus \Span(\xi _{-1}, \xi _{-2})$ 
and $\ft=\fg_{\ev}\oplus \Span(\xi _{1}, \xi _{2})$ are polarizations 
for $\lambda$.  As is easy to verify, $I^{\fg}_{\fh}(\lambda)\not\cong 
I^{\fg}_{\ft}(\lambda)$; besides, the weights of the module 
$\Ind^{\fg}_{\fh}(\lambda)$ are $\lambda$, $\lambda -x^*$, hence, 
$\lambda -(\lambda -x^*)=x^*$ but $x^*([\fg _{\od}, \fg _{\od}])\neq 
0$ contradicting the statement of Theorem 7 of \cite{K}.

The error in the proof of Theorem 7 of \cite{K} is not easy
to find: it is an incorrect induction in the proof of heading
a) on p. 80. Namely, if, in notations of \cite{K}, the subalgebra $H$
is of codimension (0, 1) then the irreducible quotients of $W$
considered as $G_{0}$-modules belong by inductive hypothesis to
one class from $L/L^{H}_{0}$, where
$$
L^{H}_{0} =\{\lambda \in \fg ^*:\lambda ([H, H])=0\}, 
$$
NOT to one class from $L/L^{G}_{0}$ as stated on p. 80, line 13 from
below.

\end{document}